\documentclass[11pt]{article}

\usepackage{amsmath} 
\usepackage{amssymb}
\usepackage{amsthm}
\usepackage[utf8]{inputenc} 
\usepackage[english]{babel} 
\usepackage{hyperref} 
\usepackage{graphicx} 
\usepackage{listings} 
\usepackage{fancyhdr}
\usepackage[margin=1.5in]{geometry}
\usepackage{xcolor,soul}
\usepackage{algcompatible}
\usepackage{algorithm}
\usepackage{algpseudocode}
\usepackage{booktabs}
\usepackage{longtable}
\usepackage{afterpage}
\usepackage{comment}
\usepackage{tikz}
\usepackage[framemethod=TikZ]{mdframed}
\usepackage{setspace}
\usepackage{mathrsfs}
\usepackage[title]{appendix}
\usepackage[T1]{fontenc}
\usepackage{authblk}
\usepackage{url}
\usepackage{enumerate}
\usepackage{framed}
\usepackage{subcaption}
\usepackage{algpseudocode, algorithm}

\makeatletter
\newcommand*{\rom}[1]{\expandafter\@slowromancap\romannumeral #1@}
\makeatother

\newtheorem{theorem}{Theorem}

\newtheorem{corollary}[theorem]{Corollary}

\newtheorem{example}[theorem]{Example}

\begin{document}


\author{Seyedmohammadhossein Hosseinian{\small $^\dagger$} and Andrew J. Schaefer{\small $^\ddagger$}} 

\affil{\small $^\dagger$University of Cincinnati, Cincinnati, OH 45221, \href{mailto:s.hosseinian@uc.edu}{\tt s.hosseinian@uc.edu}\\
$^\ddagger$Rice University, Houston, TX 77005, \href{mailto:andrew.schaefer@rice.edu}{\tt andrew.schaefer@rice.edu}}

\title{On Integer Programs with Irrational Data}
\date{December 16, 2023} 
\maketitle 
\begin{abstract}
\noindent 
An integer program (IP) with a finite number of feasible solutions may have an unbounded linear programming relaxation if it contains irrational parameters, due to implicit constraints enforced by the irrational numbers. We show that those constraints can be obtained if the irrational parameters are polynomials of roots of integers over the field of rational numbers, leading to an equivalent rational formulation. We also establish a weaker result for IPs involving the general class of algebraic irrational parameters, which extends to IPs with a particular form of transcendental numbers.

\paragraph*{Keywords:}
Integer programming, irrational parameters, philosophy of modeling.
\end{abstract}
\section{Introduction} \label{sec:intro}
Mathematical programming offers practical methods for optimization of dynamic systems, for example, through discretization and approximation of differential equations by difference equations, e.g.,~\cite{ajayi2023combination,gerdts2005solving,leyffer2022sequential,martin2006mixed,sager2011combinatorial}. Irrational numbers are ubiquitous in computational models of physical and financial systems, e.g.,~\cite{asachenkov1993disease,brunnermeier2014macroeconomic,chen1984linear,glover2016power,moller2004mixed}, which raises interest in mathematical programs in the presence of irrational parameters.

If an integer program (IP) contains only rational parameters and its linear programming relaxation (LPR) is unbounded, the IP is known to be either unbounded or infeasible~\cite{meyer1974existence}. However, if the problem involves irrational parameters, the LPR may be unbounded while the IP has a finite number of feasible solutions, and hence an optimal solution. Byrd et al.~\cite{byrd1987recognizing} provided the following example for this case:

\begin{example}  \label{ex:BGH_example}
{\em \cite{byrd1987recognizing}} 
\begin{equation*} 
\begin{array}{rl}
\max~~&  x_1\\
\text{{\em s.t.}}~~& x_3 - \sqrt{2} (x_1 - x_2) = 0,\\
& x_2 + x_4 = 1,\\
& x_1,x_2,x_3,x_4 \geq 0~~\text{{\em and integer}}. 
\end{array}
\end{equation*}
\end{example}

\medskip

\noindent
The IP has two feasible solutions $(0,0,0,1)$ and $(1,1,0,0)$ with the optimal objective value of 1, but its LPR is unbounded. This occurs because the integrality of $x_1$ and $x_2$ implies that the multiplicand of $\sqrt{2}$, i.e., $x_1 - x_2$, is rational, and if it is not zero, the left-hand side of the first constraint will be irrational while its right-hand side is rational. Thus, $x_1 - x_2 = 0$ is a logical constraint but is absent from the formulation, resulting in the unboundedness of the LPR, as seen by the ray $(1, 0,\sqrt{2}, 0)$.

Byrd et al.~\cite{byrd1987recognizing} studied IPs of the form
\begin{equation} 
\label{eq:BGH} 
\begin{array}{rl}
\max~~&  {\bf c^\top x}\\
\text{s.t.}~~& {\bf Ax = b},\\
& {\bf x \geq 0}~~\text{and integer},
\end{array}
\end{equation}
and showed that if~\eqref{eq:BGH} has only a finite number of feasible solutions but the corresponding LPR is unbounded, then its feasible region has a dimension strictly greater than that of the space spanned by the integer solutions. Irrational parameters can enforce {\em implicit} equalities among the integer (rational) solutions, the absence of which in the LPR may result in dimensionality expansion and unboundedness. While the main result of~\cite{byrd1987recognizing} establishes the existence of such implicit constraints, it does not provide a means to identify them. The elementary rationality argument employed in Example~\ref{ex:BGH_example}, to obtain the constraint $x_1 - x_2 = 0$, was due to the presence of a single irrational parameter; it fails for a constraint like 
$$
c_0 + c_1 \sqrt{2}\, x_1 + c_2 \sqrt[3]{9}\, x_2 + c_3 \pi\, x_3 = 0,
$$
where, for $i \in \{0,1,2,3\},~c_i \in \mathbb{Q}$, the set of rational numbers. To address this case, an immediate question is whether any of the irrational parameters in this constraint can be written as a linear combination of the others (over $\mathbb{Q}$) or if any linear combination of them (over $\mathbb{Q}$) may be rational.

Irrational numbers do not enjoy many favorable properties of their rational counterparts. Foremost, the set of irrational numbers is not closed under addition or multiplication. More importantly, it is nontrivial to decide if linear combinations of irrational numbers (over $\mathbb{Q}$) are irrational. For example, it is still unknown if $e+\pi$ is irrational~\cite{gowers2010princeton,klee1991old}. This makes deriving the logical (implicit) constraints for IPs involving the general class of irrational parameters extremely difficult, in the absence of further fundamental number theory results.

Here, we show that all those implicit constraints can be obtained if the irrational parameters of~\eqref{eq:BGH} are roots of integers (of arbitrary degrees) or their (multivariate) polynomials over $\mathbb{Q}$. We note that this characterization covers all polynomials of fractional powers of rational numbers (over $\mathbb{Q}$), because one can always convert the denominators to 1 within the polynomial ring. Roots of integers constitute an important and a relatively well-behaved subclass of irrational numbers, which allows for the derivation of the aforementioned result. We also show that a similar but weaker result holds for IPs involving the more general class of algebraic irrational numbers, which can be extended to IPs that involve a particular form of transcendental (non-algebraic) numbers. Our results also indicate that the new (implicit) constraints will always have rational parameters and the original constraints of~\eqref{eq:BGH} become redundant when these new constraints are added to the formulation. Therefore, the original irrational system ${\bf Ax = b}$ can be replaced by an equivalent rational system ${\bf A'x = b'}$. A direct consequence of this result is that if the LPR of such an equivalent formulation is unbounded, it implies that the original IP is either unbounded or infeasible. Equivalently, if the IP has finitely many feasible solutions, the LPR will have an optimal solution as well, generating a (finite) bound on the IP's optimal value.

\section{Preliminaries} \label{sec:pre}
Real irrational numbers consist of irrational algebraic numbers and transcendental numbers. Algebraic numbers can be expressed as the roots of polynomials with integer coefficients, while transcendental numbers cannot. All (real) transcendental numbers are irrational, because every rational number is an algebraic number of degree one. Almost all real numbers are irrational, and almost all irrational numbers are transcendental; the cardinality of the set of algebraic numbers (also rationals and integers) is $\aleph_0$, and the cardinality of the set of real numbers as well as its subset of transcendental numbers is $\aleph_1$~\cite{kaplansky2020set}. In general, deciding if a certain number is irrational, and if so, whether it is transcendental or algebraic is nontrivial. A few examples of the real numbers that are shown to be transcendental include: $e$, $\pi$, $\ln 2$, and $2^{\sqrt{2}}$~\cite{hardy1979introduction,Hermite1874,Lindemann1882}. Roots of integers (square roots, cube roots, etc.) constitute a fundamental subclass of algebraic irrational numbers. We note that the $n$-th root of an integer is either an integer or irrational~\cite{cohn1982algebra}, e.g., $\sqrt{a}$ and $\sqrt[3]{b}$ are irrational unless $a$ is a perfect square and $b$ is a perfect cube.

We briefly review a few number theoretic and algebraic concepts, as some of the results to follow are expressed in these terms. A commutative {\em ring} is a set of elements equipped with two operations (addition and multiplication), which are associative and commutative, and special elements (0 and 1) that serve as the additive and multiplicative identities. Every element $a$ of a ring has an additive inverse $-a$, i.e., $a + (-a) = 0$. A {\em field} is a commutative ring in which every nonzero element has a multiplicative inverse, i.e., $a \times \frac{1}{a} = 1$. Elements of a field can be added to, subtracted from, and multiplied by each other, and divided by any element other than 0; a field is closed under field operations. The fields of rational numbers $\mathbb{Q}$ and real numbers $\mathbb{R}$ are prominent examples. Let $K \subseteq L$ be a pair of fields. Given $\alpha \in L \backslash K$, let $K[\alpha]$ be the smallest subfield of $L$ that contains $K$ and $\alpha$. Then, by definition, $K[\alpha]$ contains all polynomials of $\alpha$ with coefficients in $K$, i.e., 
$$
c_0 + c_1 \alpha + c_2 \alpha^2 + c_3 \alpha^3 + \ldots
$$
where $c_j \in K,\forall j$. $K[\alpha]$ is called an {\em extension field} of $K$ by $\alpha$. Observe that an extension field $K[\alpha,\beta] = K[\alpha][\beta] = K[\beta][\alpha]$ contains all (multivariate) polynomials of $\alpha$ and $\beta$ with coefficients in $K$. A chain of extensions can be developed in a similar manner and may be infinitely long. In particular, the field of real numbers is an (infinite) extension field of rational numbers by the set of all irrational (real) numbers. While extension fields are theoretically well-defined, not all of them can be systematically characterized, e.g., transcendental extensions of the field of rational numbers. In the sequel, we consider the extension fields of rational numbers $\mathbb{Q}$ by a finite collection of (irrational) roots of integers as well as a finite collection of (irrational) algebraic numbers that satisfy certain properties. We refer to~\cite{morandi2012field} for a comprehensive treatment of field extensions.

\section{IPs with Irrational Parameters} \label{sec:results}
\begin{theorem} \label{th:BGH}
{\em \cite{byrd1987recognizing}} 
Let $\mathcal{P} = \{{\bf x} \in \mathbb{R}^n : {\bf A x = b}, \, {\bf x \geq 0}\}$, and $\mathcal{S}$ be the convex hull of the integer solutions in $\mathcal{P}$. If $\mathcal{P}$ is unbounded, yet contains only a finite (but nonempty) set of integer solutions, then it is of a strictly higher dimension than that of $\mathcal{S}$.  
\end{theorem}

\noindent
We note that a necessary condition for the dimensionality expansion established in Theorem~\ref{th:BGH} is that the polyhedral set defined by the constraints involving irrational parameters has an empty interior in the space of the original variables. Observe that one could not make the same argument about Example~\ref{ex:BGH_example} if the first constraint was an inequality. In fact, if a constraint defines a full-dimensional region (half space) and is written in standard form (equality), integrality/rationality on the corresponding auxiliary variable cannot be necessarily enforced. Suppose that the first constraint of Example~\ref{ex:BGH_example} were $-\sqrt{2}(x_1 - x_2) \leq 1$, written in standard form by adding $x_3$ as a slack variable. Then, any solution $(t,0,\sqrt{2}t+1,1)$, where $t$ is a nonnegative integer, would be feasible to such an IP, and both the IP and its LPR would be unbounded. The fact that Example~\ref{ex:BGH_example} has finitely many feasible solutions is because $x_3$ must be rational. It is worth mentioning that IP solution methods commonly assume integrality of auxiliary variables, which holds if the parameters are all rationals, e.g., the derivation of Gomory fractional cuts~\cite{gomory1958}.

We consider the case where the irrational parameters in the constraints of~\eqref{eq:BGH} are (multivariate) polynomials of roots of integers (of arbitrary degrees) over $\mathbb{Q}$. Equivalently, every parameter in {\bf A} or {\bf b} belongs to a finite extension of the field of rational numbers by a collection of roots of prime numbers, i.e., $\mathbb{Q}[\sqrt[\uproot{5}q_1]{p_1}, \ldots, \sqrt[\uproot{5}q_{m}]{p_m}]$, where $p_1,\ldots,p_m$ are primes and $q_1,\ldots,q_{m}$ are positive integers greater than 1. Recall that every integer greater than 1 is either a prime number or can be uniquely expressed as the product of prime numbers, by the fundamental theorem of arithmetic. We also note that, because $\mathbb{Q}[\sqrt[\uproot{5}q_1]{p_1},\ldots, \sqrt[\uproot{5}q_{m}]{p_m}]$ contains all polynomials of $\sqrt[\uproot{5}q_1]{p_1},\ldots, \sqrt[\uproot{5}q_{m}]{p_m}$ with rational coefficients, it is without loss of generality to unify the denominator of the fractional powers of the primes and associate each prime with a unique root, e.g., 
$$
1 - \frac{2}{3} (48)^{\frac{1}{6}} + \frac{1}{4} (10)^{\frac{3}{4}} = 1 - \frac{2}{3} (\sqrt[3]{2})^2 (\sqrt[6]{3}) + \frac{1}{4} (\sqrt[4]{2})^3 (\sqrt[4]{5})^3 \: \in \: \mathbb{Q}[\sqrt[3]{2},\sqrt[4]{2},\sqrt[6]{3},\sqrt[4]{5}] = \mathbb{Q}[\sqrt[12]{2},\sqrt[6]{3},\sqrt[4]{5}].
$$

\begin{theorem} \label{th:besicovitch}
{\em \cite{besicovitch1940linear}}
Let 
$$
y_1 = r_1 \, p_1,~y_2 = r_2 \,  p_2,~\ldots,~y_s = r_s \, p_s,
$$
where $p_1,p_2,\ldots,p_s$ are different primes, and $r_1,r_2,\ldots, r_s$ are positive integers not divisible by any of the these primes. If $x_1,x_2,\ldots,x_s$ are positive real roots of the equations
$$
x^{n_1} - y_1 = 0,~~x^{n_2} - y_2 = 0,~~\ldots,~~x^{n_s} - y_s = 0,
$$
and $P(x_1,x_2,\ldots,x_s)$ is a polynomial with rational coefficients of degree less than or equal to $n_1 - 1$ with respect to $x_1$, less than or equal to $n_2 - 1$ with respect to $x_2$, and so on, then $P(x_1,x_2,\ldots,x_s)$ can vanish only if its coefficients vanish. 
\end{theorem}

\begin{theorem} \label{th:main}
Consider $\mathcal{X} = \{{\bf x} \in \mathbb{R}^n : {\bf Ax = b}, \, {\bf x \geq 0}~\text{{\em and integer}}\}$, where every irrational parameter is a (multivariate) polynomial of
roots of integers with arbitrary degrees over the field of rational numbers. Then, $\mathcal{X}$ can be represented by an equivalent rational linear system.
\end{theorem}

\begin{proof}
Let $I$ denote the index set of the prime factors of all roots of integers that appear in the irrational parameters of ${\bf Ax = b}$. We consider each prime factor $p_i, \, i \in I$, to be uniquely associated with the root $\sqrt[\uproot{5}q_i]{p_i}$, where $q_i$ is the lowest common denominator of its fractional powers across all irrational numbers in {\bf A} and {\bf b}. Hence, every parameter of the system ${\bf Ax = b}$ is of the form
\begin{equation} \label{eq:form}
\sum_{k_1 = 0}^{K_1} \, \sum_{k_2 = 0}^{K_2} \, \ldots \, \sum_{k_{|I|} = 0}^{K_{|I|}} d_{k_1 k_2 \ldots k_{|I|}} \prod_{i \in I} \big ( \sqrt[\uproot{5}q_i]{p_i} \, \big )^{k_i},
\end{equation}
where $d_{k_1 k_2 \ldots k_{|I|}} \in \mathbb{Q}$, and $K_i, \, \forall i \in I$, is the largest power of $\sqrt[\uproot{5}q_i]{p_i}$ in {\bf A} and {\bf b}, modulo $q_i$. We note that the rational part of each parameter is captured by $k_1=k_2=\ldots=k_{|I|} = 0$, and for the rational parameters, this is the only term in~\eqref{eq:form} with a non-zero coefficient. Now, consider a single original equality of $\mathcal{X}$,
\begin{equation} \label{eq:const}
a_1 x_1 + \ldots + a_n x_n - b = 0,
\end{equation}
where $a_j, \forall j \in \{1,\ldots,n\}$, and $b$ are of the form~\eqref{eq:form}, and let $d_{k_1 k_2 \ldots k_{|I|}}^j$ and $d_{k_1 k_2 \ldots k_{|I|}}^b$ denote the rational coefficient in~\eqref{eq:form} for $a_j, \forall j \in \{1,\ldots,n\}$, and $b$, respectively. Observe that~\eqref{eq:const} is a polynomial of the prime roots, i.e., $P \big ( \sqrt[\uproot{5}q_1]{p_1},\ldots,\sqrt[\uproot{5}q_{|I|}]{p_{|I|}} \big )$, and the coefficient of each monomial term $\prod_{i \in I} \big ( \sqrt[\uproot{5}q_i]{p_i} \, \big )^{k_i}$, for a fixed $k_1,\ldots,k_{|I|}$, is
\begin{equation} \label{eq:coef}
\sum_{j=1}^n d_{k_1 k_2 \ldots k_{|I|}}^j x_j - d_{k_1 k_2 \ldots k_{|I|}}^b,
\end{equation}
which is rational by integrality of the ${x_j}$ variables and rationality of the $d_{k_1 k_2 \ldots k_{|I|}}^j$ and $d_{k_1 k_2 \ldots k_{|I|}}^b$ parameters. Thus, Equation~\eqref{eq:const}, as a polynomial $P \big ( \sqrt[\uproot{5}q_1]{p_1},\ldots,\sqrt[\uproot{5}q_{|I|}]{p_{|I|}} \big )$, satisfies the conditions of Theorem~\ref{th:besicovitch} (with $r_1=\ldots=r_{|I|}=1$), which implies that~\eqref{eq:coef} vanishes for every fixed $k_1,\ldots,k_{|I|}$, i.e.,
\begin{equation} \label{eq:new}
\sum_{j=1}^n d_{k_1 k_2 \ldots k_{|I|}}^j x_j = d_{k_1 k_2 \ldots k_{|I|}}^b,~\forall (k_1,\ldots,k_{|I|})\in \{0,\ldots,K_1\} \times \ldots \times \{0,\ldots,K_{|I|}\}.
\end{equation}
Recall that~\eqref{eq:new} is implied by a single equality in ${\bf A x = b}$. The first part of the proof is complete by collecting the set of constraints~\eqref{eq:new} for all defining equalities of $\mathcal{X}$ in a system ${\bf A' x = b'}$, whose parameters are all rational numbers. To show the equivalence of $\{{\bf x} \in \mathbb{R}^n : {\bf Ax = b}, \, {\bf x \geq 0}~\text{and integer}\}$ and $\{{\bf x} \in \mathbb{R}^n : {\bf A'x = b'}, \, {\bf x \geq 0}~\text{and integer}\}$, it is sufficient to note that~\eqref{eq:const} reduces to $0 = 0$ when the equalities~\eqref{eq:new} are satisfied. 
\end{proof}

\begin{corollary}
Associated with IP~\eqref{eq:BGH}, where all irrational parameters in ${\bf A x = b}$ are polynomials of roots of integers over the field of rational numbers, consider the (equivalent) IP
\begin{equation} \label{eq:ip}
\begin{array}{rl}
\max~~&  {\bf c^\top x}\\
\text{{\em s.t.}}~~& {\bf A'x = b'},\\
& {\bf x \geq 0}~~\text{{\em and integer}},
\end{array}
\end{equation}
where ${\bf A'}$ and ${\bf b'}$ are constructed according to the proof of Theorem~\ref{th:main}, involving rational parameters exclusively. Then, if~\eqref{eq:BGH} has only a finite number of feasible solutions, the LPR of~\eqref{eq:ip} gives a bound on its optimal value. Equivalently, if the LPR of~\eqref{eq:ip} is unbounded,~\eqref{eq:BGH} is either unbounded or infeasible.  
\end{corollary}

\noindent
We note that optimality or unboundedness of the LPR of~\eqref{eq:ip}, through the simplex method, is determined by the signs of the cost vector components, not their exact values. Therefore, if the vector ${\bf c}$ involves irrational parameters, the round-off-error in solving the LPR of~\eqref{eq:ip} is no worse than that of the general all-rational linear programs, as long as the irrational numbers are computed to a required precision.

Mordell~\cite{mordell1953linear} extended the linear-independence result of Besicovitch to the set of algebraic irrational numbers. He showed that, for a collection $\alpha_1,\ldots,\alpha_s$ of algebraic irrational numbers with degrees $n_1, \ldots, n_s$, respectively, the (irrational) members of $\mathbb{Q}[\alpha_1,\ldots,\alpha_s]$ are linearly independent over $\mathbb{Q}$ if 
\begin{equation} \label{eq:mordell}
\alpha_1^{\nu_1} \ldots \alpha_s^{\nu_s} \not \in \mathbb{Q},   
\end{equation}
%
for all $\nu_1,\ldots,\nu_s$ except for
$
\nu_1 \equiv 0 \, \text{(mod $n_1$)},~\ldots,~\nu_s \equiv 0 \, \text{(mod $n_s$)}.
$
Therefore, if the irrational parameters in ${\bf A x = b}$ of~\eqref{eq:BGH} are polynomials of the algebraic numbers $\alpha_1,\ldots,\alpha_s$ and when it can be verified that~\eqref{eq:mordell} holds for every $\nu_1,\ldots,\nu_s$ described above, then the coefficient of every monomial term $\prod_{i=1}^s \alpha_i^{k_i}, k_i \geq 0$, of the original system must vanish, leading to an equivalent system ${\bf A' x = b'}$ with rational parameters exclusively, akin to the proof of Theorem~\ref{th:main}.

By Lindemann–Weierstrass theorem~\cite{Baker}, a similar result holds for IPs involving a particular class of transcendental numbers. If $\alpha_1,\ldots,\alpha_s$ are algebraic numbers and linearly independent over $\mathbb{Q}$, then the exponentials $e^{\alpha_1},\ldots,e^{\alpha_s}$ are transcendental and linearly independent over the algebraic numbers. Therefore, if~\eqref{eq:BGH} contains constraints of the form 
\begin{equation} \label{eq:LW}
e^{\alpha_1} \left ( \sum_{j=1}^n a_j^1 x_j + b^1 \right ) + \ldots + e^{\alpha_s} \left ( \sum_{j=1}^n a_j^s x_j + b^s \right ) = 0,
\end{equation}
where $a_j^1, \ldots, a_j^s,~\forall j \in \{1,\ldots,n\}$, and $b^1, \ldots, b^s$ are algebraic and $\alpha_1,\ldots,\alpha_s$ are linearly-independent algebraic over $\mathbb{Q}$, then the original constraint~\eqref{eq:LW} reduces to 
$$
\sum_{j=1}^n a_j^k x_j + b^k = 0,~\forall k \in \{1,\ldots,s\},
$$
which may satisfy the Mordell's condition  and imply a system of equalities with only rational parameters. Finally, we note that Lindemann–Weierstrass theorem requires $\alpha_1,\ldots,\alpha_s$ to be algebraic and not necessarily irrational. Thus, one may always consider $\alpha_1 = 0$ in~\eqref{eq:LW} to include all the previously discussed cases.

\paragraph*{Acknowledgements.}
Partial support of Office of Naval Research grant N000142112262 is gratefully acknowledged.

\bibliographystyle{abbrv}

\end{document}